\documentclass[12pt]{article}
\usepackage{amssymb}
\usepackage{amsmath}
\usepackage{indentfirst}

\newcommand{\ba}{\begin{array}}
\newcommand{\ea}{\end{array}}
\newcommand{\bi}{\begin{itemize}}
\newcommand{\ei}{\end{itemize}}
\newcommand{\bc}{\begin{center}}
\newcommand{\ec}{\end{center}}
\newcommand{\bfr}{\begin{flushright}}
\newcommand{\efr}{\end{flushright}}
\newcommand{\f}{\frac}

\newcommand{\ds}{\displaystyle}

\newcommand{\q}{\quad}

\begin{document}

\title{A note on certain inequalities for bivariate means}
\author{J\'ozsef S\'andor\\
Babe\c{s}-Bolyai University\\
Department of Mathematics\\
Str. Kog\u{a}lniceanu nr. 1\\
400084 Cluj-Napoca, Romania\\
email: jsandor@math.ubbcluj.ro}
\date{}
\maketitle

\begin{abstract}
We obtain simple proofs of certain results from paper [1].
\end{abstract}



\section{Introduction}
Let $a,b$ be two distinct positive numbers.
The power mean of order $k$ of $a$ and $b$ is defined by
$$A_k=A_k(a,b)=\left(\ds\f{a^k+b^k}{2}\right)^{1/k},\q k\ne 0$$
and
$$A_0=\lim\limits_{k\to 0}A_k=\sqrt{ab}=G(a,b).$$

Let $A_1=A$ denote also the classical arithmetic mean of $a$ and $b$, and
$$He=He(a,b)=\ds\f{2A+G}{3}=\ds\f{a+b+\sqrt{ab}}{3}$$
the so-called Heronian mean.

In the recent paper [1] the following results have been proved:
$$A_k(a,b)>a^{1-k}I(a^k,b^k)\mbox{ for } 0<k\le 1,\ b>a
\eqno(1.1)$$
$$A_k(a,b)<I(a,b)\mbox{ for } 0<k\le \ds\f{1}{2}
\eqno(1.2)$$
$$He(a^k,b^k)<A_\beta (a^k,b^k)<\ds\f{3}{2^{1/\beta }}He(a^k,b^k)\mbox{ for }
k>0,\ \beta \ge \ds\f{2}{3}
\eqno(1.3)$$\and
$$A_k<S<2^{1/k}\cdot A_k\mbox{ for } 1\le k\le 2.
\eqno(1.4)$$

In the proofs of (1.1){(1.4) the differential calculus has been used.
Our aim will be to show that, relations (1.1)-(1.4) are easy consequences of some known results.

\section{Main results}

{\bf Lemma 2.1.}
{\it The function
$f_1(k)=\left(\ds\f{a^k+b^k}{2}\right)^{1/k}=A_k(a,b)$
is a strictly increasing function of $k$;
while $f_2(k)=(a^k+b^k)^{1/k}$
is a strictly decreasing function of $k$.
Here $k$ runs through the set of real numbers.
}

{\bf Proof.}
Through these results are essentially known in the mathematical folklore,
we shall give here a proof.

Simple computations yield:
$$k^2\ds\f{f'_1(k)}{f_1(k)}=\ds\f{x\ln x+y\ln y}{x+y}
-\ln \left(\ds\f{x+y}{2}\right),
\eqno(2.1)$$
and
$$k^2 \ds\f{f'_2(k)}{f_2(k)}=\ds\f{x\ln x+y\ln y}{x+y}-\ln(x+y),
\eqno(2.2)$$
where $x=a^k>0$, $y=b^k>0$.
Since the function $f(x)=x\ln x$ is strictly convex
(indeed:
$f''(x)=\ds\f{1}{x}>0$) by
$f\left(\ds\f{x+y}{2}\right)<\ds\f{f(x)+f(y)}{2}$,
relation (2.1) implies
$f'_1(k)>0$.
Since the function $t\to \ln t$ is strictly increasing, one has
$\ln x<\ln(x+y)$ and $\ln y<\ln(x+y)$;
so $x\ln x+y\ln y<(x+y)\ln (x+y)$,
so relation (2.2) implies that $f'_2(t)>0$.
These prove the stated monotonicity properties.

{\bf Proof of (1.1).}
By the known inequality
$I<A$ we have
$$I(a^k,b^k)<A(a^k,b^k)=\ds\f{a^k+b^k}{2}.$$

Now
$$\ds\f{a^k+b^k}{2}\le a^{k-1}\left(\ds\f{a^k+b^k}{2}\right)^{1/k}$$
is equivalent with (for $1-k>0$)
$$\left(\ds\f{a^k+b^k}{2}\right)^{1/k}>a
\mbox{ or }
a^k+b^k>2a^k,$$
which is true for $b>a$.
For $k=1$ the inequality becomes $I<A$.

{\bf Proof of (1.2).}
Since $A_k$ is strictly increasing, one has
$$A_k\le A_{1/2}=\left(\ds\f{\sqrt a+\sqrt b}{2}\right)^2=\ds\f{A+G}{2}<I,$$
by a known result (see [3]) of the author:
$$I>\ds\f{2A+G}{3}>\ds\f{A+G}{2}.
\eqno(2.3)$$

{\bf Proof of (1.3).}
By the inequality
$He<A_{2/3}$ (see [2]) one has
$$He(a^k,b^k)<A_{2/3}(a^k,b^k)\le A_\beta (a^k,b^k),$$
by the first part of Lemma 2.1.

Now,
$2^{1/\beta }A_\beta (a^k,b^k)\le 2^{3/2}(a^k,b^k)$
by the second part of Lemma 2.1, and
$A_{2/3}(a^k,b^k)<\ds\f{3}{2\sqrt 2}He(a^k,b^k)$,
by (see [2])
$$A_{2/3}<\ds\f{3}{2\sqrt 2}He.
\eqno(2.4)$$

Since $2^{3/2}=2\sqrt 2$, inequality (1.3) follows.

{\bf Proof of (1.4).}
In [2] it was proved that
$$A_2<S<\sqrt 2 A_2.
\eqno(2.5)$$

Now, by Lemma 2.1 one has, as
$k\le 2$ that $A_k\le A_2<S$ and $\sqrt 2 A_2\le 2^{1/k}A_k$.
Thus, by  (2.5), relation (1.4) follows.
We note that condition $1\le k$ is not necessary.

\end{document}